\documentclass[a4paper,11pt]{article}
\usepackage{latexsym}
\usepackage{amssymb}
\usepackage{amsfonts}
\usepackage{amsmath}
\usepackage{indentfirst}
\usepackage{graphicx}
\usepackage{epsfig}
\newtheorem{prop}{Proposition}[section]
\newtheorem{teor}{Theorem}[section]
\newtheorem{lemma}{Lemma}[section]
\newtheorem{cor}{Corollary}[section]
\newcommand{\nkinN}{n,k\in \mathbf{N}}
\newcommand{\ninN}{n\in \mathbf{N}}

\newcommand{\kinN}{k\in \mathbf{N}}
\newcommand{\cvd}{\quad $\blacksquare$\bigskip}

\date{}

\author{Luca Ferrari\thanks{Dipartimento di
Sistemi e Informatica, viale Morgagni 65, 50134 Firenze, Italy
{\tt ferrari@dsi.unifi.it}}}

\title{Some combinatorics related to central binomial coefficients: Grand-Dyck
paths, coloured noncrossing partitions and signed pattern avoiding
permutations} 

\begin{document}

\maketitle

\begin{abstract}
We give some interpretations to certain integer sequences in terms
of parameters on Grand-Dyck paths and coloured noncrossing
partitions, and we find some new bijections relating Grand-Dyck
paths and signed pattern avoiding permutations. Next we transfer a
natural distributive lattice structure on Grand-Dyck paths to
coloured noncrossing partitions and signed pattern avoiding
permutations, thus showing, in particular, that it is isomorphic
to the structure induced by the (strong) Bruhat order on a certain
set of signed pattern avoiding permutations.
\end{abstract}

\section{Introduction}

Let $\mathcal{P}$ be a set of paths in the discrete plane, having
both the starting and ending points in common. Then it is natural
to consider the partial order on $\mathcal{P}$ defined by
declaring that the path $P\in \mathcal{P}$ is less than the path
$Q\in \mathcal{P}$ when $P$ lies weakly below $Q$ (weakly meaning
that the two paths can have some points in common). This point of
view has been considered in a series of papers \cite{BBFP,BF,FP},
where some order properties of certain classical sets of lattice
paths are exploited. In particular, it is shown that the class of
Dyck paths of the same length endowed with such a partial order is
actually a distributive lattice, and the same happens for Motzkin
and Schr\"oder paths.

\bigskip

The motivation of the present work comes from \cite{BBFP}, where
it is shown that the lattices of Dyck paths are order isomorphic
to the sets of 312-avoiding permutations with the induced (strong)
Bruhat order. As a byproduct, we then have that 312-avoiding
permutations of any given length possess a distributive lattice
structure. To obtain this result, a new distributive lattice
structure is introduced and studied on noncrossing partitions.
Subsequently, similar results are proved in \cite{BF} starting
from Motzkin and Schr\"oder paths. The aim of the present work is
to consider this order structure on Grand-Dyck paths (which are,
by definition, like Dyck paths, except that they are allowed to
cross the $x$-axis), and to study some of its properties, as well
as to find some related order structure on some kind of
noncrossing partitions and pattern avoiding permutations. In
trying to accomplish this project, we will come across some
bijections and formulas which we believe to be new; in particular,
in the spirit of the last section of \cite{BBFP}, we will give
certain number sequences a combinatorial interpretation in terms
of parameters on (Grand-)Dyck paths. This is essentially the
content of section \ref{numbers}.

However, our main result is contained in section \ref{posets}, and
consists of the proof that our order structure on Grand-Dyck paths
is isomorphic to the Bruhat order on some classes of signed
pattern avoiding permutations. This generalizes the above recalled
result on Dyck paths, which in fact can be seen as a
specialization of the present one. As a byproduct, we have
determined a family of pattern avoiding signed permutations such
that the induced Bruhat order gives rise to a distributive lattice
structure (and not merely a poset structure): to the best of our
knowledge, this is the first result of this nature for signed
permutations.

\bigskip

Before starting, we recall a recursive construction of Grand-Dyck
paths based on the ECO method which will be useful in section
\ref{numbers}.

\bigskip

Let $\mathcal{GD}_n$ be the set of Grand-Dyck paths of length
$2n$, that is, by definition, the set of all lattice paths
starting at the origin $(0,0)$, ending on the $x$-axis at $(2n,0)$
and using only two kinds of steps, namely $U=(1,1)$ and
$D=(1,-1)$. Dyck paths are a special subclass of Grand-Dyck paths,
which can be obtained by adding the constraint of remaining weakly
above the $x$-axis.

\bigskip

It is possible to generate Grand-Dyck paths (according to the
semilength) using the so-called ECO method. We will not give a
description of this method here, but we refer to the very detailed
survey \cite{BDLPP}. The following construction can be found in
\cite{PPR}.

Given $P\in \mathcal{GD}_n$, we construct a set of paths of
$\mathcal{GD}_{n+1}$ as follows:
\begin{itemize}
\item[-)] if the last step of $P$ is a down step, then we insert a
peak into any point of the last descent of G or a valley into the
last point of G;

\item[-)]otherwise, we insert a valley into any point of the last
ascent of $P$ or a peak into the last point of $P$.
\end{itemize}

The succession rule $\Omega$ describing the above construction is:
\begin{equation}\label{succrule}
\Omega: \left\{ \begin{array}{lll} (2)\\
(2)\rightsquigarrow (3)(3)\\ (k)\rightsquigarrow (3)^2
(4)(5)\cdots (k)(k+1)
\end{array} \right. .
\end{equation}

\bigskip

We close this section by providing a series of notations and
definitions we will frequently need throughout all the paper.

\bigskip

An infinite lower triangular matrix $A$ is called a \emph{Riordan
array} \cite{SGWW} when its column $k$ ($k=0,1,2,\ldots$) has
generating function $d(x)(xh(x))^k$, where $d(x)$ and $h(x)$ are
formal power series with $d(0)\neq 0$.

\bigskip

We will usually denote lattice paths using capital letters, such
as $P,Q,R,\ldots$. We will also make some use of a
\emph{functional notation} for paths starting and ending on the
$x$-axis: the notation $P(k)$ stands for the ordinate of the path
$P$ having abscissa $k$.

\bigskip

In a lattice path $P$, a \emph{peak} is a sequence of two
consecutive steps, the first one being an up step and the second
one being a down step. Dually, a \emph{valley} is defined by
interchanging the role of up and down steps in the definition of a
peak.

Moreover, a \emph{descent} is a sequence of consecutive down
steps, whereas an \emph{ascent} is a sequence of consecutive up
steps.

\bigskip

For a permutation $\pi$, we use the term \emph{rise} to mean a
sequence of consecutive and increasing entries of $\pi$, whereas a
\emph{fall} is a sequence of consecutive and decreasing entries of
$\pi$.

\bigskip

We will denote with $B_n$ the \emph{hyperoctahedral group of size
$n$}, i.e. the set of all permutations of $\{ 1,2,\ldots n\}$
whose elements can be possibly signed. Signed elements will simply
be overlined. A signed element will be often interpreted as a
\emph{negative} element. In this sense, we say that the
\emph{absolute value} $|x|$ of an element $x$ is that element
without its sign. Moreover, given $\pi \in B_n$, we will denote by
$|\pi |$ the permutation of $S_n$ obtained from $\pi$ by taking
the absolute values of all its elements.

\bigskip

In every poset we will deal with, the covering relation will be
denoted $\prec$.

\bigskip

The linear order on $n$ elements, also called \emph{chain} of
cardinality $n$, will be denoted $\mathcal{C}_n$.

\bigskip

A \emph{join-irreducible} of a distributive lattice $D$ is any
element $x$ which is not the minimum of the lattice and with the
property that, if $x=u\vee v$, then $x=u$ or $x=v$.

\bigskip

The \emph{spectrum} of a distributive lattice $D$ is the poset
$\mathbf{Spec}(D)$ of the join-irreducibles of $D$.

\section{Bijections and numbers}\label{numbers}

Given a double indexed sequence $(\alpha_{n,k})_{\nkinN}$, its
\emph{coloured version} is defined to be
$(\beta_{n,k})_{\nkinN}=(2^k \alpha_{n,k})_{\nkinN}$. All the
formulas we will get in the present section can be interpreted in
the same way, namely we will provide some combinatorial
interpretation for (the row sums of) the coloured versions of a
series of (not always well known) double indexed sequences using
Grand-Dyck paths.

\bigskip

Given a Dyck path $P$, a \emph{factor} of $P$ is a minimal subpath
of $P$ which is itself a Dyck path. In figure 1 a Dyck path having
4 factors is shown.

\begin{figure}[!h]\label{granddyck}
\begin{center}
\includegraphics[scale=0.3]{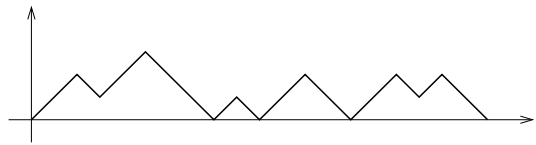}
\end{center}
\caption{A Dyck path of length 20 having 4 factors.}
\end{figure}

Now denote by $\overline{\mathcal{D}}_n$ the set of coloured Dyck
paths of length $2n$, i.e. Dyck paths whose steps can be coloured
in two different ways, say black and white. There is an obvious
bijection between $\mathcal{GD}_n$ and a special subset of
$\overline{\mathcal{D}}_n$. More precisely, we have the following,
simple proposition.

\begin{prop}\label{bijection} The set $\mathcal{GD}_n$ of Grand-Dyck paths of
length $2n$ is in bijection with the subset of
$\overline{\mathcal{D}}_n$ consisting of all coloured Dyck paths
in which steps belonging to the same factor occur with the same
colour.
\end{prop}

\emph{Proof.} For any given Grand-Dyck path, just reverse the
pieces of the path which lie below the $x$-axis and colour their
steps black (whereas the remaining steps are taken to be
white).\cvd

The subset of $\overline{\mathcal{D}}_n$ mentioned in the above
proposition will be denoted $\widetilde{\mathcal{D}}_n$ and its
elements will be called \emph{factor-bicoloured Dyck paths}. For
an example, see figure 2.

\begin{figure}[!h]\label{factbicdyck}
\begin{center}
\includegraphics[scale=0.5]{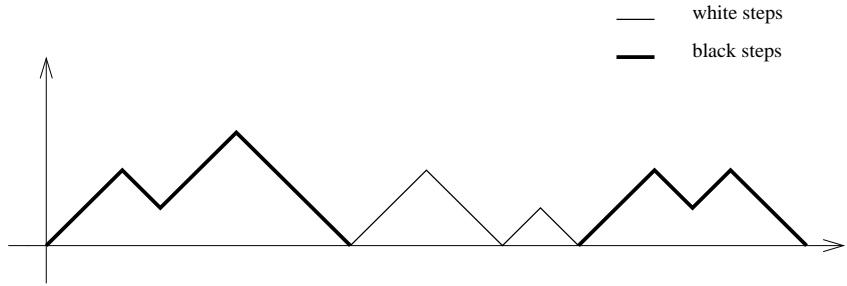}
\end{center}
\caption{A factor-bicoloured Dyck path of length 20.}
\end{figure}

\bigskip

The very easy observation expressed in the above proposition
yields the first, obvious enumerative result. Indeed, since it is
well known that Grand-Dyck path are counted by the central
binomial coefficients ${2n\choose n}$, we have the following.

\begin{prop} If $b_{n,k}=\frac{k}{2n-k}{2n-k\choose n}$ are the ballot numbers, for any $n\geq 1$, we have
\begin{equation}\label{ballot}
{2n\choose n}=\sum_{k=1}^{n}2^k b_{n,k}.
\end{equation}
\end{prop}

\emph{Proof.} To count Grand-Dyck paths of length $2n$ we can just
count the elements of $\widetilde{\mathcal{D}}_n$. Given a Dyck
path of length $2n$, the number of its factors clearly coincides
with the number of returns of the path, that is how many times the
path touches the $x$-axis except for the starting point. It is
well known (see, for example, \cite{D}) that the number of Dyck
paths of length $2n$ having precisely $k$ returns is given by the
ballot number $b_{n,k}=\frac{k}{2n-k}{2n-k\choose n}$. Since each
factor can be coloured in two different ways, the thesis
immediately follows.\cvd



The result of the above proposition is not new, and can be found,
for example, in \cite{Sl}, among the formulas for sequence A000984
(i.e. central binomial coefficients). Therefore identity
(\ref{ballot}) provides a trivial combinatorial interpretation for
the row sums of the coloured ballot numbers, that is
\begin{displaymath}
\sum_{k=1}^{n}2^k b_{n,k}=\sum_{P\in \mathcal{GD}_n}1.
\end{displaymath}

\bigskip

A definitely more interesting result can be obtained by
generalizing a result described in the last section of
\cite{BBFP}. To this aim, we need first of all to introduce a
bijection between Grand-Dyck paths and a special class of set
partitions.

If we denote by $\Pi_n$ the set of partitions of $\{ 1,2,\ldots
,n\}$, given $\pi =B_1 |B_2 |\cdots |B_h \in \Pi_n$, we will
always represent it in such a way that (i) the elements inside
each block $B_i$ are listed in decreasing order and (ii) the
blocks are listed in increasing order of their maxima. This will
be called the \emph{standard representation, or standard form} of
$\pi$. Thus, for instance, the partition $\{ 4\} |\{ 5,3,1\} |\{
6,2\}$ in $\Pi_6$ is represented here in its standard form. To
improve readability, in the sequel we will also delete all the
parentheses and commas, so that the above partition will be
written $4|531|62$. Each partition can be factored as follows.
Given $\pi =B_1 |B_2 |\cdots |B_h \in \Pi_n$ represented in its
standard form, we say that \emph{$\pi$ has $k$ components} when
its (linearly ordered) set of blocks can be partitioned into $k$
nonempty intervals with the property that the union of the blocks
inside the same interval is an interval of $\{ 1,2,\ldots ,n\}$,
and $k$ is maximum with respect to this property. Each of the
above intervals of blocks will be called a \emph{component} of
$\pi$. For instance, the partition $2|43|651|8|97$ in $\Pi_9$ has
2 components, which are $2|43|651$ and $8|97$. We will denote
$t_{n,k}$ the number of partitions of an $n$-set having $k$
components. Set partitions having a single component are also
called \emph{atomic partitions} (see \cite{BZ}).

\begin{prop} The infinite triangular array $(t_{n,k})_{\nkinN}$ is
the Riordan array $(p(x),p(x))$, where $p(x)$ is the generating
function of atomic partitions. In particular, $t_{n,k}$ is the
coefficient of $x^n$ in $x^k p(x)^k$.
\end{prop}

\emph{Proof.}\quad Observe that a partition of an $n$-set having
$k$ components can be uniquely recovered by the subpartition
constituted by its first $k-1$ components and the atomic partition
isomorphic to its $k$-th component. This argument can be
translated into the following recurrence relation:
\begin{displaymath}
t_{n,k}=\sum_{h=1}^{n}t_{n-h,k-1}p_h ,
\end{displaymath}
where $p_n$ denotes the number of atomic partitions of an $n$-set.
If $C_k (x)$ is the generating function of the $k$-th column of
the array $T=(t_{n,k})_{\nkinN}$, the above recurrence becomes:
\begin{displaymath}
C_k (x)=xp(x)C_{k-1}(x),
\end{displaymath}
where $p(x)=\sum_{n>1}p_n x^n$ is the generating function of
atomic partitions. Iterating we then get:
\begin{displaymath}
C_k (x)=(xp(x))^k ,
\end{displaymath}
which is precisely our thesis.\cvd

The sequence of atomic partitions is also recorded in \cite{Sl}
(it is sequence A074664), and its generating function is
$p(x)=1-\frac{1}{B(x)}$, where $B(x)$ is the (ordinary) generating
function of Bell numbers (see, for example, \cite{Kl}). The
infinite matrix $T$ is in \cite{Sl} too (sequence A127743), but
the combinatorial interpretation given here is different.

\bigskip

\emph{Remark.}\quad Observe that the generating function $p(x)$
can also be determined using a species-theoretic argument. Indeed,
the fact that any nonempty partition can be decomposed into the
(possibly empty) partition constituted by all but the last of its
components and the atomic partition constituted by the last
component alone means that the species of nonempty partitions is
obtained as the Hadamard product of the species of partitions and
the species of atomic partitions, and so we get the generating
function relation:
\begin{displaymath}
B(x)-1=B(x)\cdot p(x),
\end{displaymath}
whence the equality $p(x)=1-\frac{1}{B(x)}$ follows.

\bigskip

To generalize the last result of \cite{BBFP}, which gives a
combinatorial interpretation of Bell numbers in terms of natural
parameters on Dyck paths, we now need to introduce the notion of
\emph{component-bicoloured partition}. As the name itself
suggests, a \emph{component-bicoloured partition} is a set
partitions whose components can be coloured using two different
colours, say black and white. The total number of
component-bicoloured partitions of an $n$-set is clearly given by
$\sum_{k=1}^{n}2^k t_{n,k}$. This is sequence A059279 in
\cite{Sl}, but also in this case the present interpretation is not
recorded. We now find a further combinatorial interpretation of
this sequence in terms of natural parameters on Grand-Dyck paths.

A \emph{bicoloured Dyck word} is a Dyck word (i.e. a word on the
alphabet $\{ U,D\}$ such that, interpreting each $U$ as an up step
and each $D$ as a down step, the resulting path is a Dyck path)
whose letters can be coloured either black or white. We call
\emph{factor-bicoloured Dyck word} any bicoloured Dyck word
corresponding to a factor-bicoloured Dyck path.


Define a \emph{bicoloured Bell matching} of a factor-bicoloured
Dyck word to be any Bell matching of the associated Dyck word
(i.e. the word obtained by ``forgetting colours"). Following
\cite{BBFP}, a Bell matching of a Dyck word $\omega$ is a matching
between the $U$'s and the $D$'s of $\omega$ such that
\begin{enumerate}
\item for any set of consecutive $D$'s, the leftmost $D$ is
matched with the adjacent $U$ on its left;
\item every other $D$ is matched with a $U$ on its left, in such a
way that there are no crossings among the arcs originated from a
set of consecutive $D$'s.
\end{enumerate}

Observe that, in a bicoloured Bell matching, if a $U$ and a $D$
are matched, then they have the same colour.

Adapting the argument given for proposition 6.1 in \cite{BBFP},
the reader can now easily prove the following.

\begin{prop} There is a bijection between bicoloured Bell
matchings of factor-bicoloured Dyck words of length $2n$ and
component-bicoloured partitions of an $n$-set.
\end{prop}

Given two bicoloured Bell matchings, we say that they are
\emph{equivalent} when they are bicoloured Bell matchings of the
same factor-bicoloured Dyck word. Since it is clear that, for each
fixed factor-bicoloured Dyck word there is exactly one bicoloured
Bell matching without crossings among its arcs, as an immediate
consequence of the last proposition we have that in the
equivalence classes of the above defined equivalence relation
there is precisely one bicoloured Bell matching corresponding to a
component-bicoloured noncrossing partition. It is convenient to
record this fact in a proposition.

\begin{prop}\label{bij1} There is a bijection between $\mathcal{GD}_n$ (or,
which is the same, $\widetilde{\mathcal{D}}_n$) and the set
$\widetilde{NC}(n)$ of component-bicoloured noncrossing partitions
of an $n$-set.
\end{prop}

Now we are ready to state the main result of this section, which
provides a combinatorial interpretation for the coloured version
of the row sums of the sequence $(t_{n,k})_{\nkinN}$ using
Grand-Dyck paths. In the following theorem, the word ``positive"
means ``above the $x$-axis", and ``negative" stands for ``below
the $x$-axis". Moreover, the \emph{(absolute) height} of a peak
(or a valley) is given by the (absolute value of) the ordinate if
its vertex.

\begin{teor} Given a Grand-Dyck path $P$, let $\mathcal{A}$ be the
set of positive peaks and negative valleys of $P$ and
$\mathcal{B}$ the set of non-positive peaks, non-negative valleys
and returns on the $x$-axis of $P$. If we linearly order the
elements of $\mathcal{A}$ and $\mathcal{B}$ using their abscissas,
and denote by $p_1 ,p_2 ,\ldots ,p_h$ the absolute heights of the
elements of $\mathcal{A}$ and by $v_1 ,v_2 ,\ldots ,v_h$ the
absolute heights of the elements of $\mathcal{B}$ (with the
convention $v_h =0$), then we have
\begin{equation}\label{main}
\sum_{k=1}^{n}2^k t_{n,k}=\sum_{P\in
\mathcal{GD}_n}\prod_{i=1}^{h}{p_i -1\choose v_i}.
\end{equation}
\end{teor}

\emph{Proof.}\quad First of all, we observe that
$|\mathcal{B}|=|\mathcal{A}|$, since each element of $\mathcal{A}$
is followed by precisely one element of $\mathcal{B}$. Therefore
the statement of the theorem is proved to be consistent.

Now let $\pi$ be a component-bicoloured noncrossing partition of
an $n$-set; proposition \ref{bij1} implies that $\pi$ is uniquely
associated with a Grand-Dyck path $P$, and so with a
factor-bicoloured Dyck path $\widetilde{P}$ of length $2n$. If we
denote by $Q$ the Dyck path obtained from $\widetilde{P}$ by
``forgetting" colours, we observe that the $p_i$'s are the heights
of the peaks of $Q$, whereas the $v_i$'s are the heights of the
valleys (except for $v_h$, which is the height of the last return
of $Q$, and so $v_h =0$). Therefore, recalling the definition of
Bell matching and the last theorem of \cite{BBFP}, if $|[\pi ]|$
is the equivalence class of $\pi$, we obtain:
\begin{displaymath}
|[\pi ]|=\prod_{i=1}^{h}{p_i -1\choose v_i}.
\end{displaymath}

Summing over all Grand-Dyck paths, we then get formula
(\ref{main}), and the theorem is proved.\cvd

To conclude the part of this section devoted to partitions, we
observe that the notion of component-bicoloured partition is a
specialization of the well-known notion of \emph{bicoloured
partition}, i.e. partition with bicoloured blocks (see, for
instance, \cite{LLPP}). As it is obvious, such partitions are
enumerated by the sequence $\sum_{k=1}^{n}2^k S_{n,k}$, where the
$S_{n,k}$'s are the Stirling numbers of the second kind; this is
sequence A001861 in \cite{Sl}. Also in this case there is a
combinatorial interpretation of this numbers in terms of
Grand-Dyck paths. The key result is the following proposition,
which refines the formula found in \cite{BBFP} to express Bell
numbers in terms of parameters on Dyck paths and whose proof
(which can be carried out by suitably generalizing the argument of
\cite{BBFP}) is left to the reader.

\begin{prop} If $\mathcal{D}_n (k)$ denotes the set of Dyck paths
of length $2n$ having exactly $k$ peaks, it is
\begin{displaymath}
S_{n,k}=\sum_{P\in \mathcal{D}_n (k)}\prod_{i=1}^{k}{p_i -1\choose
v_i},
\end{displaymath}
where $p_i$ and $v_i$ are the heights of the peaks and the valleys
of $P$, respectively (with $v_k =0$ by convention).
\end{prop}


As an immediate consequence, we have the following alternative
interpretation for the row sums of the coloured Stirling numbers
of the second kind.

\begin{cor} With the same notations as in the above proposition, we
have:
\begin{displaymath}
\sum_{k=1}^{n}2^k S_{n,k}=\sum_{k=1}^{n}2^k \sum_{P\in
\mathcal{D}_n (k)}\prod_{i=1}^{k}{p_i -1\choose v_i}.
\end{displaymath}
\end{cor}

\bigskip

The second part of the present section is devoted to the
description of some new bijections between Grand-Dyck paths and
signed pattern avoiding permutations. More precisely, we propose
here two bijections: the former has been found with the help of
the ECO method, whereas the latter will be useful in the next
section to define an interesting distributive lattice structure.

\bigskip

The first bijection involves the classes of signed pattern
avoiding permutations $B_n (21,\overline{21})$, where $B_n$ is the
hyperoctahedral group on $n$ elements. It is known \cite{Si} that
$|B_n (21,\overline{21})|={2n\choose n}$. Moreover, combining some
results in \cite{Si} and in \cite{Re}, an explicit bijection
between $\mathcal{GD}_n$ and $B_n (21,\overline{21})$ can be
described. However, our bijection is different from the one so
obtained. Before starting, observe that a permutation in $B_n
(21,\overline{21})$ is a shuffle of the signed and unsigned
elements, each of which are ordered increasingly by absolute
value. Thus, given $\pi \in B_n (21,\overline{21})$, we can
consider the two elements $a$ and $b$ such that $|a|$ and $|b|$
are the maximum of the absolute values of the signed and of the
unsigned elements, respectively. It is clear that $|a|=n$ or
$|b|=n$. We call \emph{quasi maximum} of $\pi$ the one between $a$
and $b$ whose absolute value is different from $n$.

To describe our first bijection, we represent permutations by a
graphical device used, for instance, in \cite{BFP}. We represent
the elements of the permutations as dots placed on horizontal
lines, in such a way that elements with greater absolute values
lie on higher lines. It is an extremely natural representation, so
we deem it is not necessary to give a more formal definition (see
figure 3 for an example).

\begin{figure}[!h]\label{staff}
\begin{center}
\includegraphics[scale=0.5]{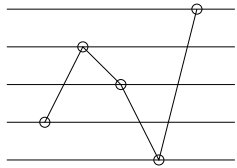}
\end{center}
\caption{A graphical representation of the permutation 24315}
\end{figure}

Let $\pi \in B_n (21,\overline{21})$. Starting from $\pi$ we will
construct a set of permutations belonging to
$B_{n+1}(21,\overline{21})$, by adding a new element in the last
position of $\pi$ and then suitably renaming some of the elements
of $\pi$. From a graphical point of view, we simply add a new
horizontal line in the representation of $\pi$, and we place on
such a line the new element. In performing this operation, we have
to take care that the resulting permutations still avoids the two
patterns 21 and $\overline{21}$. To make sure that this happens,
we can distinguish two cases.
\begin{enumerate}
\item Suppose that the quasi maximum $\overline{a}$ of $\pi$ is
signed. In this case, we can add at the end of $\pi$ any signed
element having absolute value greater than $a$, as well as the
unsigned element $n+1$. Figure 4 describes how this construction
works.

\begin{figure}[!h]\label{ecoperm1}
\begin{center}
\includegraphics[scale=0.4]{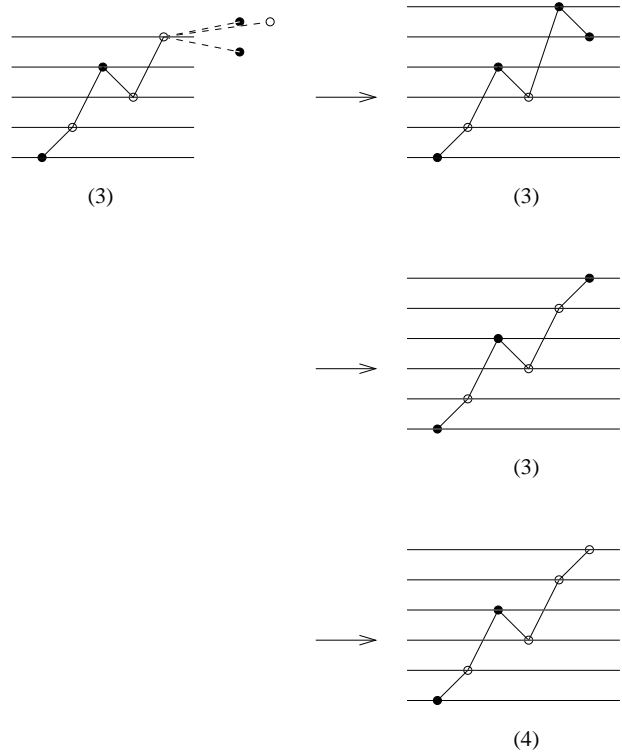}
\end{center}
\caption{Our ECO construction performed on
$\overline{1}2\overline{4}35$. Here signed elements are
represented using black bullets.}
\end{figure}

\item On the other hand, if the quasi maximum $a$ of $\pi$ is
unsigned, we are allowed to add at the end of $\pi$ any unsigned
element greater than $a$, as well as the signed element
$\overline{n+1}$ (see figure 5).

\begin{figure}[!h]\label{ecoperm2}
\begin{center}
\includegraphics[scale=0.4]{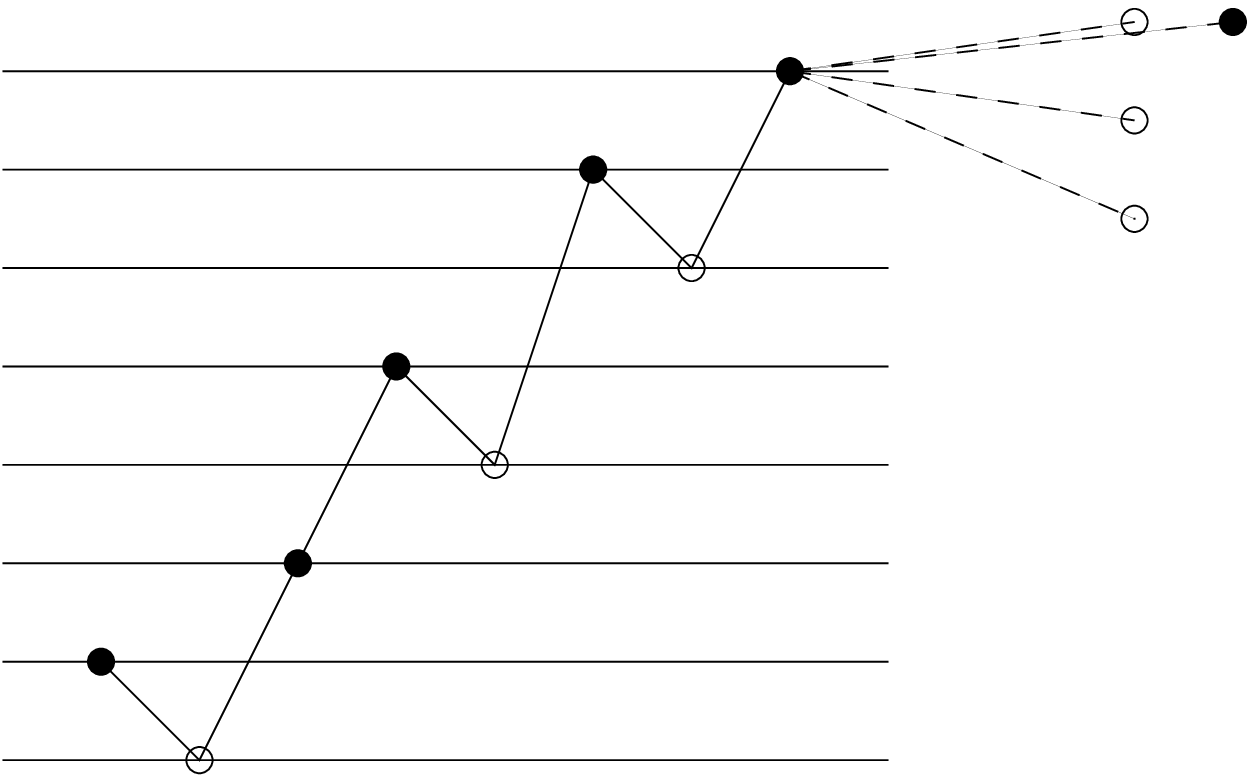}
\end{center}
\caption{Our ECO construction performed on
$\overline{2}1\overline{35}4\overline{7}6\overline{8}$.}
\end{figure}

\end{enumerate}

In the first case, if $\overline{a}=\overline{n+2-k}$, then $\pi$
produces $k$ sons, whose quasi maximums are easily seen to be
$\overline{n+2-k},\overline{n+3-k},\ldots ,\overline{n},n$, and so
the numbers of their sons are, respectively, $k+1,k,\ldots 4,3,3$.
Similarly, in the second case, we have the same statement as
above, with signed elements replaced by unsigned ones and vice
versa. Thus, also in this case the numbers of sons are the same as
above.

Therefore, we observe that we have an ECO construction for signed
permutations avoiding 21 and $\overline{21}$ which is isomorphic
to the ECO construction for Grand-Dyck paths recalled in the
introduction, and described by the succession rule $\Omega$ given
in (\ref{succrule}). Such an isomorphism defines a bijection
between $\mathcal{GD}_n$ and $B_n (21,\overline{21})$. Moreover,
the underlying ECO construction allows to easily translate many
natural statistics on Grand-Dyck paths into specific parameters
defined on permutations. To give just a glimpse of this fact, we
propose a couple of examples.

\bigskip

Given a signed permutation $\pi =\pi_1 \cdots \pi_n \in B_n$, we
define the following two sets:
\begin{eqnarray*}
A(\pi )=\{ i\leq n\; |\; \textnormal{$\pi_i$ is unsigned and
$\pi_1 \cdots \pi_{i-1}$ has a signed maximum}\}\cup \nonumber
\\ \{ i\leq n\; |\;
\textnormal{$\pi_i$ is signed and $\pi_1 \cdots \pi_{i}$ has an
unsigned maximum}\} ;
\\ B(\pi )=\{ i\leq n\; |\; \textnormal{the two greatest elements of $\pi_1
\cdots \pi_i$ have different signs}\} .
\end{eqnarray*}

By convention, we assume that $1\in A(\pi )$ if and only if
$\pi_1$ is unsigned and that $1\in B(\pi )$ for any permutation
$\pi$.

\begin{prop} The statistic ``number of peaks" on Grand-Dyck paths
corresponds to the A-statistic on $B_n (21,\overline{21})$.
\end{prop}

\emph{Proof.}\quad Observe that, in the ECO construction of
Grand-Dyck paths proposed in the introduction (and encoded by the
rule in (\ref{succrule})), the sons of a Grand-Dyck path $P$
having $h$ peaks can have either $h$ peaks or $h+1$ peaks. More
precisely:
\begin{itemize}
\item[\emph{$path_1$})] if $P$ ends with a sequence of down steps,
then all the sons of $P$ have $h+1$ peaks but the one obtained by
adding a valley at the end of $P$ and the one obtained by adding a
peak at the beginning of the last descent of $P$; if $P$ has label
$(k)$, the two sons of $P$ having $h$ peaks are then labelled
$(3)$ and $(k+1)$;

\item[\emph{$path_2$})]if $P$ ends with a sequence of up steps,
then all the sons of $P$ have $h+1$ peaks but the one obtained by
adding a valley at the beginning of the last ascent of $P$ (if $P$
has label $(k)$, then the label of such a son is $(k+1)$).
\end{itemize}

Transferring the above considerations on permutations by means of
our bijection, we have the following two cases:
\begin{itemize}
\item[\emph{$perm_1$})] if $\pi \in B_n (21,\overline{21})$ has an
unsigned maximum, then, in the above described ECO construction of
$B_n (21,\overline{21})$, the two sons of $\pi$ having labels
$(3)$ and $(k+1)$ corresponding to the paths mentioned in
\emph{$path_1$}) are avoided by adding a signed element $a$ whose
absolute value is not greater than the absolute values of all the
elements of $\pi$; this means that $\pi a$ has an unsigned
maximum;

\item[\emph{$perm_2$})] if $\pi \in B_n (21,\overline{21})$ has a
signed maximum, then the son of $\pi$ labelled $(k+1)$ is avoided
by adding any signed element.
\end{itemize}

The above considerations immediately implies that, if $\pi$
corresponds to $P$ in our bijection, then $h=|A(\pi )|$, which is
the thesis.\cvd

\begin{prop} The statistic ``number of returns" on Grand-Dyck paths
corresponds to the B-statistic on $B_n (21,\overline{21})$.
\end{prop}

\emph{Proof (sketch).}\quad The arguments to be used here are
completely analogous to those of the previous proposition. We just
observe that, in the ECO construction of Grand-Dyck paths, a new
return is produced whenever either a valley or a peak is appended
at the end of the path. This translates on permutations into the
addition to the right of $\pi$ of an element $a$ whose sign is
different from the sign of the maximum of $\pi$ and such that the
maximum of $\pi$ and $a$ are the two greatest elements of $\pi
a$.\cvd

In passing through, we notice that it is possible to define
another (presumably new) bijection between $\mathcal{GD}_n$ and
$B_n (21,\overline{21})$. Also in this case, we start by
considering the usual ECO construction of Grand-Dyck paths and
then we translate it into permutations avoiding the two patterns
$21$ and $\overline{21}$. Without going into details, the idea is
to generate permutations by adding a new maximum (instead of
adding the rightmost element). This can be represented by using a
graphical device similar to the one above: just replace horizontal
lines with vertical lines. In figure 6 an example of how this
construction works is shown. We entirely leave to the interest
reader the accomplishment of all the details of this alternative
approach, as well as the task of translating some statistics on
Grand-Dyck paths (such as the number of peaks and the number of
returns considered above) into permutations.

\begin{figure}[!h]\label{ecoperm4}
\begin{center}
\includegraphics[scale=0.4]{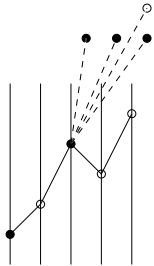}
\end{center}
\caption{An alternative ECO construction performed on
$\overline{1}2\overline{4}35$.}
\end{figure}

\bigskip

Our second bijection involves a different type of pattern avoiding
permutations, namely permutations which avoid the four patterns
$312,\overline{312},2\overline{1},\overline{2}1$. The classes of
pattern-avoiding permutations $B_n
(312,\overline{312},2\overline{1},\overline{2}1)$ are counted by
the central binomial coefficients ${2n\choose n}$, and this can be
easily proved by exhibiting a completely trivial bijection with
component-bicoloured noncrossing partitions.


\begin{prop}\label{barremoving}(Bar-removing bijection). There is a bijection between
$\widetilde{NC}(n)$ and $B_n
(312,\overline{312},2\overline{1},\overline{2}1)$.
\end{prop}

\emph{Proof.}\quad Taken $\pi \in \widetilde{NC}(n)$, written as
usual in its standard form, delete the vertical bars so to obtain
a permutation belonging to $B_n$ (still denoted by $\pi$). The
presence of a pattern $2\overline{1}$ or of a pattern
$\overline{2}1$ in $\pi$ would imply that, in the associated
partition, the elements of such a pattern should belong to two
different components and the greatest of them should belong to a
block with a lesser index. But this is impossible in our standard
representation of partitions. Moreover, the fact that $|\pi |$ is
noncrossing implies that every signed version of the pattern 312
cannot appear in the associated signed permutation. Finally, it is
immediate that, avoiding any signed version of 312 and the two
patterns $\overline{2}1$ and $2\overline{1}$ is equivalent to
avoiding $312,\overline{312},\overline{2}1$ and $2\overline{1}$.

To prove that this is actually a bijection, it is sufficient to
observe that, given $\pi \in B_n
(312,\overline{312},2\overline{1},\overline{2}1)$, the associated
partition can be uniquely recovered by inserting a vertical bar
between the elements of each rise of $|\pi |$ \cvd

Now, using propositions \ref{bij1} and \ref{barremoving}, we get
the following corollary.

\begin{cor}There is a bijection between $\mathcal{GD}_n$ and $B_n
(312,\overline{312},2\overline{1},\overline{2}1)$.
\end{cor}


This last bijection is a sort of signed analog of a well known
bijection between Dyck paths and 312-avoiding permutations, which
can be found for example in \cite{BK,Kr}. It also has the
remarkable feature of translating many natural statistics on paths
into natural statistics on permutations. For instance, the number
of unsigned (resp. signed) left-to-right maxima in a permutation
of $B_n (312,\overline{312},2\overline{1},\overline{2}1)$ is equal
to the number of positive peaks (resp. negative valley) of the
associated Grand-Dyck path.

\section{Posets}\label{posets}

The set $\mathcal{GD}_n$ of Grand-Dyck paths of length $2n$ can be
naturally ordered by declaring $P\leq Q$ whenever $P(k)\leq Q(k)$,
for all $k$. This means that the path $P$ lies weakly below $Q$
(see the example in figure 7). This very natural partial order is
easily seen to yield a distributive lattice structure, in which
the join and meet of two paths are taken coordinatewise. The
minimum and maximum of these lattices will be denoted $\mathbf{0}$
and $\mathbf{1}$, respectively. See figure 8 for the Hasse diagram
of $\mathcal{GD}_3$. The resulting lattice structures have already
appeared in the literature (see \cite{NF}), but they have been
considered on different combinatorial objects. It is immediate the
following fact, whose easy proof is left to the reader.

\begin{figure}[!h]\label{ordergd}
\begin{center}
\includegraphics[scale=0.4]{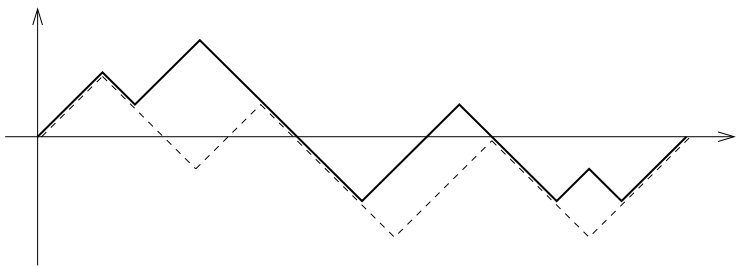}
\end{center}
\caption{Two comparable Grand-Dyck paths of length 20.}
\end{figure}

\begin{figure}[!h]\label{latticegd}
\begin{center}
\includegraphics[scale=0.3]{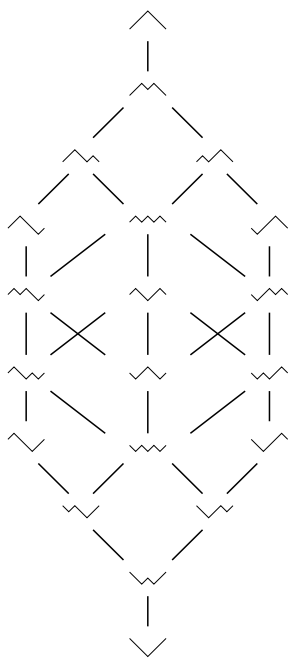}
\end{center}
\caption{The lattice $\mathcal{GD}_3$.}
\end{figure}

\begin{prop} The lattice $\mathcal{GD}_n$ of Grand-Dyck paths of length
$2n$ is isomorphic to the Young lattice of integer partitions
which lie inside the $n\times n$ square.
\end{prop}

Young lattices of partitions have been intensively investigated
since many years, see for instance \cite{St} for an interesting
study on the unimodality properties of such lattices. Of course,
as an immediate consequence of the above proposition, we have that
Grand-Dyck lattices are rank-unimodal.

However, even if the abstract lattice structure we are considering
is not new, we claim that the study of order and lattice
properties arising from the special representation in terms of
Grand-Dyck paths is worth being carried out.

\bigskip

Our first result concerns the shape of join-irreducible elements.

\begin{prop}
A Grand-Dyck path $P\in \mathcal{GD}_n$ is join-irreducible if and
only if it has precisely one peak. Therefore,
$\mathbf{Spec}(\mathcal{GD}_n)\simeq \mathcal{C}_n^2$.
\end{prop}

\emph{Proof (sketch).}\quad The covering relation on Grand-Dyck
paths works as follows: $Q$ covers $P$ if and only if $Q$ can be
obtained from $P$ by changing a valley into a peak. Then the first
part of the thesis immediately follows. As far as the second part
is concerned, just observe that, thanks to the last proposition,
join-irreducibles corresponds to integer partitions of rectangular
shape.\cvd

We have already mentioned that the lattices $\mathcal{GD}_n$ are
rank-unimodal. The rank function $r_n$ of $\mathcal{GD}_n$ is
clearly related to the area function. More precisely, we have the
following proposition.

\begin{prop}
If $P\in \mathcal{GD}_n$, then
\begin{displaymath}
r_n (P)=\frac{A(P)+n^2}{2},
\end{displaymath}
where $A(P)$ denotes the area (with sign) of the region included
between the path and the $x$-axis.
\end{prop}

\emph{Proof.}\quad By induction, suppose that $P\prec Q$ in
$\mathcal{GD}_n$ and that $r_n (P)=\frac{A(P)+n^2}{2}$. Since $Q$
is obtained from $P$ by simply reversing a valley, we have that
$A(Q)=A(P)+2$, and so $r_n
(Q)=\frac{A(Q)+n^2}{2}=\frac{A(P)+2+n^2}{2}=\frac{A(P)+n^2}{2}+1=r_n
(P)$. Since $r_n
(\mathbf{0})=\frac{A(\mathbf{0})+n^2}{2}=\frac{-n^2 +n^2}{2}=0$,
the proof is completed.\cvd

Our next goal will be to translate the above described lattice
structure on partitions. To accomplish this task we make use of
the bijection between the set $\widetilde{\mathcal{D}}_n$ of
factor-bicoloured Dyck paths of length $n$ and
component-bicoloured noncrossing partitions stated in proposition
\ref{bij1}. For the sake of simplicity, from now on we will refer
to black (resp., white) items (steps, factors,...) as to coloured
(resp., noncoloured) items.

We start by observing that the Grand-Dyck lattice structure can be
read off on factor-bicoloured Dyck paths as follows: if $P,Q\in
\widetilde{\mathcal{D}}_n$, we say that $P\leq Q$ when, for every
$\kinN$, one of the following holds:
\begin{enumerate}
\item if $P(k)$ and $Q(k)$ both belong to coloured factors, then
$P(k)\geq Q(k)$;

\item if $P(k)$ and $Q(k)$ both belong to noncoloured factors,
then $P(k)\leq Q(k)$;

\item $P(k)$ belongs to a coloured factor and $Q(k)$ belongs to a
noncoloured one.
\end{enumerate}

Now, if we transport this lattice structure along the above
recalled bijection, we obtain a lattice structure on
component-bicoloured noncrossing partitions. In the sequel, we
denote these lattices of partitions $\widetilde{NC}(n)$ (where $n$
is the size of the ground set, of course).

Our next results concerns the description of the order relation on
$\widetilde{NC}(n)$.

\bigskip

Given $\pi \in \widetilde{NC}(n)$, define the \emph{max-vector} of
$\pi$ to be the vector $max(\pi )$ whose $i$-th component equals
the maximum of the first $i$ element of $|\pi |$, when $\pi$ is
written in its standard form; moreover, each component of $max(\pi
)$ appears coloured when it is coloured in $\pi$. Therefore, for
instance, taken $\pi
=2|431|\overline{6}|\overline{7}|\overline{985}\in
\widetilde{NC}(9)$, we have $max(\pi
)=(2,4,4,4,\overline{6},\overline{7},\overline{9},\overline{9},\overline{9})$.

For any given $\ninN$, denote by $M(n)$ the set of max-vectors of
$\widetilde{NC}(n)$, that is $M(n)=\{ v=(v_1 ,\ldots ,v_n )\; |\;
\exists \pi \in \widetilde{NC}(n):v=max(\pi )\}$. It is not
difficult to see that $M(n)$ consists of all vectors with $n$
bicoloured components having increasing absolute values and such
that, for any $i<n$, $|v_i |\geq i$ and, if $v_i$ and $v_{i+1}$
have different colours, then $v_i =i$.

Now define on $M(n)$ a partial order as follows. Given two
max-vectors $v=(v_1 ,\ldots ,v_n )$ and $w=(w_1 ,\ldots ,w_n )$,
we say that $v\leq w$ when, for every $i\leq n$ either
\begin{itemize}
\item[(i)] $v_i$ and $w_i$ are both coloured and $v_i \geq w_i$,
or

\item[(ii)] $v_i$ and $w_i$ are both noncoloured and $v_i \leq
w_i$, or

\item[(iii)] $v_i$ is coloured and $w_i$ is noncoloured.
\end{itemize}

Clearly, the covering relation of the poset $[M(n);\leq ]$ can be
described by saying that precisely one of the above situations
$(i),(ii)$ or $(iii)$ holds for a specific $i$, whereas all the
other components are equal, and in each of the three cases we have
respectively:
\begin{itemize}
\item[(i)]$v_i =w_i +1$,

\item[(ii)]$w_i =v_i +1$,

\item[(iii)]$|v_i |=|w_i |$.
\end{itemize}

Using max-vectors it is now possible to characterize the covering
relation of $\widetilde{NC}(n)$.

\begin{prop}\label{maxvectors} Given $\pi ,\rho \in \widetilde{NC}(n)$, it is $\pi
\prec \rho$ if and only if $max(\pi )\prec max(\rho )$.
\end{prop}

\emph{Proof.}\quad Saying that $\pi \prec \rho$ in
$\widetilde{NC}(n)$ means that, if we consider the two associated
factor-bicoloured Dyck paths $P=P(\pi )$ and $R=R(\rho )$, they
differ precisely in two steps, namely either:
\begin{itemize}
\item[(i)] a coloured peak of $P$ is changed into a coloured
valley of $R$, or

\item[(ii)] a noncoloured valley of $P$ is changed into a
noncoloured peak of $R$, or

\item[(iii)] a coloured peak of $P$ is changed into a noncoloured
peak of $R$ (in this case the two peaks necessarily lies on the
$x$-axis).
\end{itemize}

Now observe that, given $v=(v_1 ,\ldots ,v_n )=max(\pi )\in M(n)$,
the position $d_i (\pi )$ of the $i$-th down step in the
associated factor-bicoloured Dyck path is given by $v_i +i$, and
of course the same happens for the max-vector $w=(w_1 ,\ldots ,w_n
)=max(\rho )\in M(n)$ (i.e., $d_i (\rho )=w_i +i$). Therefore, in
the above three cases, we have:
\begin{itemize}
\item[(i)] $P$ and $R$ coincide, except for a pair of adjacent
steps, which is $\overline{UD}$ in $P$ and $\overline{DU}$ in $R$.
If the down step involved is the $i$-th, then we have $w_i =d_i
(\rho )-i=d_i (\pi )-1-i=v_i -1$.

\item[(ii)] $P$ and $R$ coincide, except for a pair of adjacent
steps, which is $DU$ in $P$ and $UD$ in $R$. If the down step
involved is the $i$-th, then we have $w_i =d_i (\rho )-i=d_i (\pi
)+1-i=v_i +1$.

\item[(iii)] $P$ and $R$ coincide, except for a pair of adjacent
steps, which is $\overline{UD}$ in $P$ and $UD$ in $R$, and such a
peak lies on the $x$-axis. In this last case, if the down step
involved is the $i$-th, then we have that $v_i$ and $w_i$ have the
same absolute value, but $v_i$ is coloured whereas $w_i$ is
noncoloured.
\end{itemize}

Thus, the fact that $\pi \prec \rho$ in $\widetilde{NC}(n)$ is
equivalent to the fact that $max(\pi )\prec max(\rho )$ in $M(n)$,
which is our thesis.\cvd

\begin{cor} Given $\pi ,\rho \in \widetilde{NC}(n)$, it is $\pi
\leq \rho$ if and only if $max(\pi )\leq max(\rho )$.
\end{cor}

\bigskip

Our last goal is to transfer the above order on signed pattern
avoiding permutations. This can be done by simply applying the
bar-removing bijection of proposition \ref{barremoving}, thanks to
which we obtain a partial order on $B_n
(312,\overline{312},2\overline{1},\overline{2}1)$. The following
proposition gives a characterization of the associated covering
relation.

We will use the term \emph{signed} (respectively, \emph{unsigned})
\emph{inversion} of $\pi$ to mean a pair $(\pi_i ,\pi_j )$, with
$i<j$, $\pi_i =\overline{a}$ (resp., $\pi_i =a$), $\pi_j
=\overline{b}$ (resp., $\pi_j =b$) and $a>b$. An analogous
definition is given for the term \emph{signed (unsigned)
noninversion}.

\begin{prop}\label{permcov} Let $\pi ,\rho \in B_n
(312,\overline{312},2\overline{1},\overline{2}1)$. Then $\pi \prec
\rho$ if and only if $\rho$ is obtained from $\pi$ by either:
\begin{itemize}
\item[(i)] interchanging the elements of a signed inversion, or

\item[(ii)] interchanging the elements of an unsigned
noninversion, or

\item[(iii)]changing a signed element belonging to a rise of $|\pi
|$ into an unsigned one.
\end{itemize}
\end{prop}

\emph{Proof.}\quad The condition $\pi \prec \rho$ in $B_n
(312,\overline{312},2\overline{1},\overline{2}1)$ can be
translated on factor-bicoloured paths by saying that exactly one
of the three conditions listed in the proof of proposition
\ref{maxvectors} holds. Using the results of \cite{BBFP}, it is
not difficult to conclude that the first two conditions
corresponds, on permutations, to conditions $(i)$ and $(ii)$ of
the present proposition. Concerning the third condition, if, in a
path $P$, a coloured peak at height 0 is changed into a
noncoloured one, then, on the associated permutation $\pi$, we
have a signed element belonging to a rise of $|\pi |$ which is
changed into an unsigned one, and the proof is completed.\cvd

The distributive lattice structure on $B_n
(312,\overline{312},2\overline{1},\overline{2}1)$ obtained via the
bar-removing bijection turns out to have a very interesting
alternative combinatorial description. Before giving it
explicitly, we determine a formula to compute its rank function.
Given a permutation $\pi \in B_n
(312,\overline{312},2\overline{1},\overline{2}1)$, we denote by
$inv(\pi )$ the number of \emph{unsigned} inversions of $\pi$, by
$ninv(|\pi |)$ the number of non-inversions of $|\pi |$ and by $\#
(\pi )$ the number of \emph{unsigned} elements of $\pi$. Thus, for
instance, given $\pi =2431\overline{675}8\overline{9}\in
\widetilde{NC}(9)$, we have $inv(\pi )=4$, $ninv(|\pi |)=30$ and
$\# (\pi )=5$.

\begin{prop} If $\pi \in B_n
(312,\overline{312},2\overline{1},\overline{2}1)$, denoting by
$r_n (\pi )$ its rank, we have:
\begin{displaymath}
r_n (\pi )=ninv(|\pi |)+2inv(\pi )+\# (\pi ).
\end{displaymath}
\end{prop}

\emph{Proof.}\quad The minimum $\mathbf{0}$ of the lattice $B_n
(312,\overline{312},2\overline{1},\overline{2}1)$ is the
permutation $\overline{n(n-1)\cdots 21}$, and it is clear that
$ninv(|\overline{n(n-1)\cdots 21|})+2inv(\overline{n(n-1)\cdots
21})+\# (\overline{n(n-1)\cdots 21})=0=r_n (\mathbf{0})$.

Using an induction argument, suppose that, in $B_n
(312,\overline{312},2\overline{1},\overline{2}1)$, we have $\pi
\prec \rho$. Thanks to proposition \ref{permcov}, we have three
possible cases, which we deal with using the same numeration as in
the statement of such a proposition.
\begin{itemize}
\item[(i)] In this case, in the associated coloured Dyck path
there is a coloured peak which is changed into a coloured valley.
Referring to \cite{BBFP}, we can then say that, in the permutation
$\pi$, a new signed noninversion is produced, and so:
\begin{displaymath}
inv(\rho )=inv(\pi ),\quad ninv(|\rho |)=ninv(|\pi |)+1,\quad \#
(\rho )=\# (\pi ),
\end{displaymath}
whence
\begin{eqnarray*}
ninv(|\rho |)+2inv(\rho )+\# (\rho )&=&ninv(|\pi |)+2inv(\pi )+\#
(\pi )+1
\\&=&r_n (\pi )+1=r_n (\rho ).
\end{eqnarray*}
\item[(ii)] This situation corresponds to changing a valley in a
peak in a noncoloured factor of the associated coloured Dyck path.
But is it known from \cite{BBFP} that this produces one more
unsigned inversion in $\rho$ (leaving unchanged all the signed
elements), and so
\begin{displaymath}
inv(\rho )=inv(\pi )+1,\quad ninv(|\rho |)=ninv(|\pi |)-1,\quad \#
(\rho )=\# (\pi ),
\end{displaymath}
and an analogous computation as above immediately yields
\begin{displaymath}
ninv(|\rho |)+2inv(\rho )+\# (\rho )=r_n (\rho ).
\end{displaymath}
\item[(iii)] Finally, in this last case the difference between
$\pi$ and $\rho$ consists of the fact that $\rho$ contains one
more unsigned element, whence (using analogous arguments as those
for the preceding two cases) the thesis follows.\cvd
\end{itemize}

Our last result, which we deem is the main one of the present
paper, is the determination of an isomorphism between our lattices
$B_n (312,\overline{312},2\overline{1},\overline{2}1)$ and an
important and well known poset structure on permutations. To this
aim, we need to introduce a few notations and definitions.

Call $S_{\pm n}$ the set of permutation of the set $\{
\overline{n},\overline{n-1},\ldots
,\overline{2},\overline{1},1,2,\ldots ,n-1,n\}$, linearly ordered
as indicated. As already remarked, this corresponds to setting
$\overline{k}=-k$ and then considering the usual linear order. In
what follows, we will often tacitly make the above identification,
but we keep on writing $\overline{k}$ instead of $-k$ in order to
gain a better readability. Given $\pi$ in the hyperoctahedral
group $B_n$, denote by $\hat{\pi }$ the permutation of $S_{\pm n}$
defined by $\hat{\pi }(i)=\overline{\pi (n+1-i)}$ and $\hat{\pi
}(\overline{i})=\overline{\hat{\pi}(i)}$, for every $i\in \{
1,\ldots ,n\}$. Thus, for instance, if $\pi
=3\overline{2}5\overline{41}\in B_5$, then
$\hat{\pi}=3\overline{2}5\overline{41}14\overline{5}2\overline{3}\in
S_{\pm 5}$. Observe, in particular, that, given $\hat{\pi }$
expressed, as usual, in one-line notation, $\pi$ is obtained,
again in one-line notation, by taking the first half of the
elements of $\hat{\pi}$. Moreover, let $\hat{S}_{\pm
n}[312,\overline{312},2\overline{1},\overline{2}1]=\{ \hat{\pi }\;
|\; \pi \in B_n
(312,\overline{312},2\overline{1},\overline{2}1)\}$ (here we use
square brackets in order to avoid confusion with the notation for
pattern avoidance).

\begin{teor} The poset $B_n
(312,\overline{312},2\overline{1},\overline{2}1)$ is isomorphic to
$\hat{S}_{\pm n}[312,\overline{312},2\overline{1},\overline{2}1]$
endowed with the (strong) Bruhat order.
\end{teor}

\emph{Proof.}\quad Suppose first that $\pi \prec \rho$ in $B_n
(312,\overline{312},2\overline{1},\overline{2}1)$. From
proposition \ref{permcov}, we have three cases to analyze.
\begin{itemize}
\item[(i)] Suppose that $\rho$ is obtained from $\pi$ by
interchanging the elements of a signed inversion. To fix the
notations, we can denote by $i,j\leq n$ the two elements such that
$i<j$, $\pi (i)=\rho (j)$ and $\pi (j)=\rho (i)$ are both signed
and $\pi (i)>\pi (j)$ (and so $\rho (i)<\rho (j)$). Then, from the
definition of $\hat{\pi}$ and $\hat{\rho}$, it follows that
$\hat{\pi}(\overline{n+1-i})=\hat{\rho}(\overline{n+1-j})$ and
$\hat{\pi}(\overline{n+1-j})=\hat{\rho}(\overline{n+1-i})$ are
both signed and
$\hat{\pi}(\overline{n+1-i})>\hat{\pi}(\overline{n+1-j})$ (and so
$\hat{\rho}(\overline{n+1-i})<\hat{\rho}(\overline{n+1-j})$).
Since $n+1-i>n+1-j$, this means that $\hat{\rho}$ possesses at
least one more inversion than $\hat{\pi}$. An analogous argument
on the elements $\hat{\pi}(n+1-i)=\hat{\rho}(n+1-j)$ and
$\hat{\pi}(n+1-j)=\hat{\rho}(n+1-i)$ shows that $\hat{\rho}$ has
one further inversion more than $\hat{\pi}$. It is then easy to
realize that these are the only inversions of $\hat{\rho}$ which
are not also in $\hat{\pi}$, and so we can conclude that
$\hat{\rho}$ has two more inversions than $\hat{\pi}$ in the
Bruhat order of $\hat{S}_{\pm
n}[312,\overline{312},2\overline{1},\overline{2}1]$, which is
enough to say that $\hat{\pi}\leq \hat{\rho}$ in such a Bruhat
poset.

\item[(ii)] An analogous argument can be developed when $\rho$ is
obtained from $\pi$ by interchanging the elements of an unsigned
noninversion. Also in this case, following the same lines, it is
possible to show that $\hat{\rho}$ has two more inversions than
$\hat{\pi}$ in the Bruhat order of $\hat{S}_{\pm
n}[312,\overline{312},2\overline{1},\overline{2}1]$.

\item[(iii)] If $\rho$ is obtained from $\pi$ by changing a signed
element belonging to a rise into an unsigned one, then $\hat{\pi}$
and $\hat{\rho}$ coincide except for two elements; more precisely,
there exists a positive $i\leq n$ such that
$\hat{\pi}(\overline{i})=\hat{\rho}(i)$ is signed (and so
$\hat{\pi}(i)=\hat{\rho}(\overline{i})$ is unsigned). Thus, the
pair $(\hat{\pi}(\overline{i}),\hat{\pi}(i))$ is a noninversion in
$\hat{\pi}$, whereas $(\hat{\rho}(\overline{i}),\hat{\rho}(i))$ is
an inversion in $\hat{\rho}$. From this we deduce that
$\hat{\rho}$ has one more inversion than $\hat{\pi}$.
\end{itemize}

The above arguments allows us to conclude that, if $\pi \prec
\rho$ in $B_n (312,\overline{312},2\overline{1},\overline{2}1)$,
then $\hat{\pi}\leq \hat{\rho}$ in the Bruhat order of
$\hat{S}_{\pm n}[312,\overline{312},2\overline{1},\overline{2}1]$.
As an obvious consequence, we have that, if $\pi \leq \rho$ in
$B_n (312,\overline{312},2\overline{1},\overline{2}1)$, then
$\hat{\pi}\leq \hat{\rho}$ in the Bruhat order of $\hat{S}_{\pm
n}[312,\overline{312},2\overline{1},\overline{2}1]$, which is the
first part of the theorem.

Vice versa, suppose that $\hat{\pi}\prec \hat{\rho}$ in
$\hat{S}_{\pm n}[312,\overline{312},2\overline{1},\overline{2}1]$
with the induced Bruhat order. This means that $\hat{\rho }$ is
obtained from $\hat{\pi }$ by performing as little inversions as
possible (and, of course, remaining inside the class $\hat{S}_{\pm
n}[312,\overline{312},2\overline{1},\overline{2}1]$). We can
distinguish some cases.

\begin{itemize}
\item Given $i>j>0$, if $\hat{\pi }(\overline{i})<\hat{\pi
}(\overline{j})$ and both $\hat{\pi }(\overline{i}),\hat{\pi
}(\overline{j})$ are unsigned, then we can exchange $\hat{\pi
}(\overline{i})$ and $\hat{\pi }(\overline{j})$, provided that
there is no $k$, $i>k>j$, such that $\hat{\pi
}(\overline{i})<\hat{\pi }(\overline{k})<\hat{\pi
}(\overline{j})$. However, in this case, we also need to exchange
the two elements $\hat{\pi }(i)$ and $\hat{\pi }(j)$ in order to
remain inside $\hat{S}_{\pm
n}[312,\overline{312},2\overline{1},\overline{2}1]$. Thus, we have
obtained a permutation $\hat{\rho }$ having two more inversions
than $\hat{\pi }$. Translating all this on $B_n
(312,\overline{312},2\overline{1},\overline{2}1)$, we have that
$\rho$ is obtained from $\pi$ by interchanging the elements of an
unsigned noninversion. Of course, we have, by symmetry, exactly
the same situation if we start by considering $0<i<j$ such that
$\hat{\pi }(i)<\hat{\pi }(j)$ and $\hat{\pi }(i),\hat{\pi }(j)$
both signed.

\item With a completely analogous argument, we can prove that a
permutation $\hat{\rho }$ which covers $\hat{\pi }$ in the induced
Bruhat order of $\hat{S}_{\pm
n}[312,\overline{312},2\overline{1},\overline{2}1]$ can be
obtained by interchanging $\hat{\pi }(i)$ and $\hat{\pi }(j)$,
when $0<i<j$, $\hat{\pi }(i)$ and $\hat{\pi }(j)$ are both
unsigned and $\hat{\pi }(i)<\hat{\pi }(j)$ (or, equivalently, when
$i>j>0$, and $\hat{\pi }(\overline{i})<\hat{\pi }(\overline{j})$).
In $B_n (312,\overline{312},2\overline{1},\overline{2}1)$, this
means that $\rho$ is obtained from $\pi$ by interchanging the
elements of a signed inversion.

\item If $\hat{\pi }(i)$ and $\hat{\pi }(j)$ have different signs,
for $i,j>0$, we cannot exchange $\hat{\pi }(i)$ and $\hat{\pi
}(j)$, since a $\overline{2}1$ pattern would arise. And the same
would happen for $\hat{\pi }(\overline{i})$ and $\hat{\pi
}(\overline{j})$.

\item The only case which does not fit into one of the above is
when $\hat{\rho }$ is obtained from $\hat{\pi }$ by interchanging
two elements $\hat{\pi }(\overline{i})$ and $\hat{\pi }(j)$, with
$i,j>0$. In this case, it is easy to see that, if $\hat{\pi
}(\overline{i})$ and $\hat{\pi }(j)$ had different absolute
values, then, in $\rho \in B_n
(312,\overline{312},2\overline{1},\overline{2}1)$ we would have
two elements having the same absolute value, which is clearly not
allowed. The only possibility we have to perform an interchange is
to have $i=j$ (i.e., to interchange two elements of the kind
$a,\overline{a}$). In this case, $\rho$ is obtained from $\pi$ by
changing a signed element into an unsigned one, and, in order that
the inversion in $\hat{S}_{\pm
n}[312,\overline{312},2\overline{1},\overline{2}1]$ is minimal, it
is necessary that, for every $0<k<i$, $\hat{\pi }(k)>\hat{\pi
}(i)>0$. This means that, in $B_n$, $\pi (n+i-1)$ is signed and
belongs to a rise of $|\pi |$.
\end{itemize}

Now, putting things together, thanks to proposition \ref{permcov},
we have shown that, if $\hat{\pi }\prec \hat{\rho }$ in
$\hat{S}_{\pm n}(312,\overline{312},2\overline{1},\overline{2}1)$,
then $\pi \prec \rho$ in $B_n
(312,\overline{312},2\overline{1},\overline{2}1)$, which is enough
to conclude.\cvd

This last result is a ``signed generalization" of the fact (proved
in \cite{BBFP}) that the lattices of Dyck paths are isomorphic to
the lattices of $312-$avoiding permutations under the Bruhat
order. Indeed, consider the hyperoctahedral group $B_n$ endowed
with the Bruhat order, as it is defined, for instance, in
\cite{BB}. Using our language, it can be described as follows.
Given $\pi',\rho'\in B_n$, consider the permutations
$\hat{\pi},\hat{\rho}\in S_{\pm n}$ defined by the juxtaposition
of $\pi$ and $\pi'$ and of $\rho$ and $\rho'$, respectively, where
$\pi$ (resp., $\rho$) is defined by reversing and changing all the
signs of $\pi'$ (resp., $\rho$). For instance, if
$\pi'=14\overline{5}2\overline{3}\in B_5$, then $\pi
=3\overline{2}5\overline{41}$ and
$\hat{\pi}=3\overline{2}5\overline{41}14\overline{5}2\overline{3}\in
S_{\pm 5}$. Then $\pi'\leq \rho'$ in the Bruhat order of $B_n$ if
and only if $\hat{\pi}\leq \hat{\rho}$ in the Bruhat order of the
symmetric group $S_{\pm n}$. Therefore, as a consequence of the
last theorem, we get our final results, which states that $B_n
(312,\overline{312},2\overline{1},\overline{2}1)$ is isomorphic to
a set of signed pattern avoiding permutation under the Bruhat
order.

\begin{cor} The poset $B_n
(312,\overline{312},2\overline{1},\overline{2}1)$ is isomorphic to
$B_n (213,\overline{213},1\overline{2},\overline{1}2)$ endowed
with the Bruhat order.
\end{cor}

\emph{Proof.}\quad Just observe that, if $\pi$ and $\pi'$ are
related as above, then $\pi$ avoids a pattern $\sigma$ if and only
if $\pi'$ avoids the pattern obtained by reversing $\sigma$ and
changing the signs of all its elements.\cvd



\begin{figure}[!h]\label{latticebruhat}
\begin{center}
\includegraphics[scale=0.4]{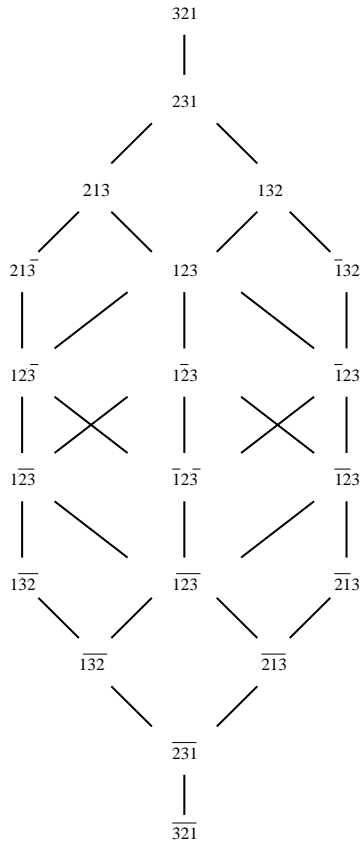}
\end{center}
\caption{The lattice $B_3
(312,\overline{312},2\overline{1},\overline{2}1)$}
\end{figure}

\end{document}